\documentclass[12pt]{article}
\usepackage{amsmath}
\usepackage{latexsym}
\usepackage{amssymb}
\newcommand{\wb}{\overline}
\newcommand{\at}{\symbol{'100}}

\newcommand{\R}{{\mathbb R}}
\newcommand{\sub}{\subseteq}
\DeclareMathOperator{\Supp}{supp}
\begin{document}
$\;$\\[-30mm]
\begin{center}
{\Large\bf Products of locally compact spaces are\\[.3mm]
{\boldmath$k_\R$}-spaces}\\[3mm]
{\bf Helge Gl\"{o}ckner and Niku Masbough}\vspace{-1mm}
\end{center}
\begin{abstract}
\hspace*{-6.3mm}A theorem by N. Noble from 1970 asserts that every product
of completely regular, locally pseudo-compact $k_\R$-spaces is
a $k_\R$-space. As a consequence,
all direct products of locally compact Hausdorff spaces
are $k_\R$-spaces. We provide a streamlined proof
for this fact.
\end{abstract}
Stimulated by \cite{Nob}, we prove the following theorem:\\[3mm]
{\bf Theorem.}
\emph{Let $(X_\alpha)_{\alpha\in A}$ be a family
of topological spaces $X_\alpha$ such that each $x\in X_\alpha$
has a neighbourhood basis consisting of closed, quasi-compact
$x$-neighbourhoods.
Let $Y$ be a topological space such that each $y\in Y$
has a neighbourhood basis consisting of closed $y$-neighbourhoods.\footnote{These conditions
are satisfied, e.g., if each $X_\alpha$ is locally compact and $Y$ is a regular
(but not necessarily completely regular) topological space. See also~\cite[p.\,187, lines 11-13]{Nob}.}
Endow $X:=\prod_{\alpha\in A}X_\alpha$ with the product
topology and let $f\colon X\to Y$ be a map such that $f|_K$ is continuous
for each quasi-compact
subset $K\sub X$. Then~$f$ is continuous.}\\[3mm]
According to~\cite{Nob}, a topological space~$X$ is called
a \emph{$k_\R$-space} if~$X$ is completely regular
and functions $f$ from~$X$ to~$\R$ (and hence also functions $f$ from $X$
to a completely regular topological space~$Y$) are continuous if and only if
$f|_K$ is continuous for each compact subset~$K\sub X$.
Since products of completely regular spaces are completely regular and every
locally compact space is completely regular, the preceding theorem entails
the following fact, which is a special case of \cite[Theorem~5.6\,(ii)]{Nob}:\\[2.3mm]
{\bf Corollary.}
\emph{For every family $(X_\alpha)_{\alpha\in A}$
of locally compact topological spaces, the product topology makes
$X:=\prod_{\alpha\in A} X_\alpha$
a $k_\R$-space.}\\[2.3mm]
Noble leaves the proof of \cite[Theorem~5.6\,(ii)]{Nob}
to the reader, and asserts that it follows by an ``obvious adaptation of the proof of 5.3''. We did not find this
assertion
convincing in the general
case of \emph{loc.\,cit.} Yet, the proof of our theorem is an adaptation of
the one of \cite[Theorem~5.3]{Nob}. It varies a less general version in~\cite{Mas}.
The corollary has been used in recent research, \cite{Nik}.\\[2.3mm]
Our terminology is as in \cite{Bou}. Notably, locally compact spaces are Hausdorff.\\[2.3mm]
{\bf Proof of the theorem.}
Let $x=(x_\alpha)_{\alpha\in A}\in X$.
To see that $f$ (as described in the theorem) is continuous at~$x$,
let us show that $f^{-1}(U)$ is a neighbourhood of~$x$ in~$X$
for each neighbourhood $U$ of~$f(x)$ in~$Y$.
By hypothesis, there exists an open neighbourhood $V$ of~$f(x)$ in~$Y$
whose closure $\wb{V}$ is contained in~$U$.
For each $\alpha\in A$, let $K_\alpha$ be a quasi-compact
closed neighbourhood of~$x_\alpha$ in~$X_\alpha$.
Since
$K:=\prod_{\alpha\in A} K_\alpha$
is quasi-compact by Tychonoff's Theorem, $f|_K$ is continuous by hypothesis.
Hence $(f|_K)^{-1}(V)$ is a neighbourhood of~$x$ in~$K$.
As~$K$ is endowed with the product topology, we find neighbourhoods
$L_\alpha$ of~$x_\alpha$ in~$K_\alpha$ such that
$\prod_{\alpha\in A}L_\alpha\sub f^{-1}(V)$
and $N:=\{\alpha\in A\colon L_\alpha\not=K_\alpha\}$
is finite. After shrinking~$L_\alpha$ for $\alpha\in N$,
we may assume that each $L_\alpha$ is closed in~$X_\alpha$ and quasi-compact.
Replacing $K_\alpha$ with $L_\alpha$ for all $\alpha$,
we may assume that
\begin{equation}\label{fK}
f(K)\sub V.
\end{equation}
For $y=(y_\alpha)_{\alpha\in A}\in X$,
let us write
$\Supp_K(y):=\{\alpha\in A\colon y_\alpha\not\in K_\alpha\}$.
For each finite subset $F\sub A$, let $\Sigma_F(K)$ be the set of all
$y=(y_\alpha)_{\alpha\in A}\in X$ such that
$\Supp_K(y)$ is finite and
$\Supp_K(y)\cap F=\emptyset$,
i.e.,
$y_\alpha\in K_\alpha$ for all $\alpha\in F$.\\[2.4mm]
{\bf Claim\,({\boldmath$*$}):} We claim that there exists a finite subset $F\sub A$ such that
\begin{equation}\label{propclai}
\Sigma_F(K)\subseteq f^{-1}(\overline{V}).
\end{equation}
If this is true, then\vspace{-9mm}
\begin{equation}\label{the-nbhd}
\prod_{\alpha\in F}K_\alpha \times \prod_{\alpha\in A\setminus F} X_\alpha
\sub f^{-1}(U),\vspace{-2mm}
\end{equation}
whence $f^{-1}(U)$ is a neighbourhood of~$x$.
In fact, if $y=(y_\alpha)_{\alpha\in A}$ is an element of the left-hand side
of~(\ref{the-nbhd}), then
\[
y^G_\alpha:=\left\{
\begin{array}{cl}
y_\alpha &\mbox{if $\,\alpha\in G$,}\\
x_\alpha & \mbox{else}
\end{array}
\right.
\]
defines an element $y^G:=(y^G_\alpha)_{\alpha\in A}\in\Sigma_F(K)$
for each $G$ in the set of\linebreak
finite subsets of~$A$ containing~$F$
(which is directed via inclusion).
Since
$C:=\prod_{\alpha\in A} (K_\alpha\cup\{y_\alpha\})$
is quasi-compact, $f|_C$ is continuous.
As the net $y^G$ converges to~$y$ in~$C$ with respect
to the product topology, we deduce that
$f(y^G)\to f(y)$.
Since $f(y^G)\in \overline{V}$ for all~$G$ (see (\ref{propclai})),
we deduce that $f(y)\in\wb{V}\sub U$, establishing~(\ref{the-nbhd}).\\[2.3mm]
To prove Claim\,($*$), let us suppose it was false
and derive a contradiction. If the claim was false,
we could obtain a sequence
$(x^n)_{n\in{\mathbb N}}$
of elements $x^n=(x^n_\alpha)_{\alpha\in A}\in X$
such that the following holds for all $n\in{\mathbb N}$:\vspace{.3mm}

(a)
$\Supp_K(x^n)$ is finite;

(b) $\Supp_K(x^m)\cap \Supp_K(x^\ell)=\emptyset$ for all $m,\ell\leq n$ such that $m\not=\ell$;

(c)
$f(x^n)\in Y \setminus \overline{V}$.\vspace{.3mm}

\noindent
If this is true, we define
$K_\alpha':=K_\alpha\cup\{x^n_\alpha\colon n\in{\mathbb N}\}$
and note that $K_\alpha'\setminus K_\alpha$ is either empty or a singleton,
by~(b).
Hence $K_\alpha'$ is quasi-compact (like $K_\alpha$).
As $K':=\prod_{\alpha\in A}K_\alpha'$
is quasi-compact and $(x^n)_{n\in{\mathbb N}}$
is a sequence in~$K'$, we see that $(x^n)_{n\in{\mathbb N}}$
has a convergent subnet $(x^{n(j)})_{j\in J}$
for some directed set~$J$. Let~$z=(z_\alpha)_{\alpha\in A}\in K'$ be a limit
of the subnet.
Given~$\alpha$, condition~(b) implies that $x^n_\alpha\in K_\alpha$ for all sufficiently
large~$n$,
whence also $x^{n(j)}_\alpha\in K_\alpha$ eventually.
Since~$K_\alpha$ is closed in~$X_\alpha$ and $x^{n(j)}_\alpha\to z_\alpha$, we deduce that
$z_\alpha\in K_\alpha$
for all~$\alpha$ and hence~$z\in K$.
But
\begin{equation}\label{agailit}
f(x^{n(j)})\in Y\setminus \overline{V}\subseteq Y \setminus V.
\end{equation}
Since $f|_{K'}$ is continuous, we have
$f(x^{n(j)})\to f(z)$. As the right-hand side $Y\setminus V$ of~(\ref{agailit})
is closed, we deduce that
$f(z)\in Y \setminus V$.
But $f(z)\in f(K)\subseteq V$ by (\ref{fK}), contradiction.\\[2.3mm]
It only remains to construct $x^1,x^2,\ldots$,
which we achieve by recursion.
As we suppose that Claim\,($*$) is false, $\Sigma_\emptyset(K)$
is not a subset of $f^{-1}(\overline{V})$;
we therefore find an element
$x^1\in\Sigma_\emptyset(K)$ such that $f(x^1)\not\in\overline{V}$.
If $x^1,\ldots, x^n$ satisfying (a)--(c) have been found, then
$F:=\Supp_K(x^1)\cup\cdots\cup\Supp_K(x^n)$
is a finite set. As we assume that Claim\,($*$) is false,
there exists $x^{n+1}\in\Sigma_F(K)$ such that $f(x^{n+1})\not\in\overline{V}$.
Since $\Supp_K(x^m)\subseteq F$ for $m\in\{1,\ldots, n\}$ but
$\Supp_K(x^{n+1})\subseteq A\setminus F$ as $x^{n+1}\in\Sigma_F(K)$,
we have $\Supp_K(x^m)\cap\Supp_K(x^{n+1})=\emptyset$.
As condition~(b) already holds for~$n$, the preceding shows that it also holds
for $n+1$ in place of~$n$.
Note that $\Supp_K(x^{n+1})$ is finite since $x^{n+1}\in\Sigma_F(K)$.
This completes the recursive construction.\vspace{-5.3mm}
{\small{\bf Helge  Gl\"{o}ckner}, Institut f\"{u}r Mathematik, Universit\"at Paderborn,\\
Warburger Str.\ 100, 33098 Paderborn, Germany.
Email: {\tt  glockner\at{}math.upb.de}\\[2.3mm]
{\bf Niku Masbough},
Institut f\"{u}r Mathematik, Universit\"at Paderborn,\\
Warburger Str.\ 100, 33098 Paderborn, Germany.
Email: {\tt nikumasbough\at{}yahoo.com}}\vfill
\end{document}